\documentclass{amsart}
\makeatletter
\def\LaTeX{\leavevmode L\raise.42ex
   \hbox{\kern-.3em\size{\sf@size}{0pt}\selectfont A}\kern-.15em\TeX}
\makeatother

\newcommand{\BibTeX}{{\rm B\kern-.05em{\sc
i\kern-.025emb}\kern-.08em\TeX}}

\theoremstyle{definition}

\numberwithin{equation}{section}

\begin{document}

\title{Reconstruction of Paley-Wiener functions on the Heisenberg group}

\author{Isaac Pesenson}
\address{Department of Mathematics, Temple University,
Philadelphia, PA 19122} \email{pesenson@math.temple.edu}


\maketitle

\begin{abstract}

  Let $M$ be a Riemmanian manifold with bounded geometry.
  We consider a generalization of Paley-Wiener functions and Lagrangian
  splines on $M$.
 An analog of the Paley-Wiener theorem is given. We also show that every
Paley-Wiener function on a manifold is uniquely
determined by its values on some discrete sets of points.

 The main result of the paper is a generalization of the Whittaker-Shannon
formula for reconstruction of a Paley-Wiener function from its values on a
discrete set. It is shown that every Paley- Wiener function on $M$ is a limit
of some linear combinations of fundamental solutions of the powers of the
Laplace-Beltrami operator.

 The result is new even in the one-dimentional case.

\end{abstract}

{ \bf 1. } We introduce an appropriate generalization of Paley-Wiener functions. Our goal is to show that the reconstruction of such functions is possible as long as the distance between points from a discrete subgroup is small enough. The reconstruction formula involves the notion of a spline. In classical case this approach was used by Schoenberg [14].

 Remark that consideration in the present paper is subelliptic in the sense that central role belongs to a sertain subelliptic operator. The case of corresponding elliptic theory on manifolds was considerd by author in [12], [13].

 The Heisenberg group $H_{m}$ is the Lie group whose underlying manifold is direct product of $R$ and $C^{m}$ and composition is given by the formula
$$(t,z)(t^{'},z^{'})=(t+t^{'}+2Im zz^{'},z+z^{'})$$ where $t,t^{'}\in R,z,z^{'}\in C^{m}$.
Dilations on $H_{m}$ are given by
$$\delta_{s}(t.z)=(s^{2}t,sz)$$ and homogeneous norm by
$$|(t,z)|=(t^{2}+|z|^{4})^{1/4}.$$

 The number $Q=2m+2$ is called the homogeneous dimention of $H_{m}$ and in analysis on $H_{m}$ plays the same role as usual dimention $d$ in analysis on $R^{d}$.

 To introduce Fourier transform on $H_{m}$ we consider irreducible unitary representations of $H_{m}$ in the Bargmann space which consist of holomorphic functions $F$ on $C^{m}$ such that
$$\|F\|^{2}=(2\lambda/\pi)^{m}\int_{C^{m}}|F(w)|^{2} exp (-2\lambda|w|^{2})dw$$ is finite. The monomials
$$F_{\alpha,\lambda}(w)=((2\lambda)^{1/2}w)^{\alpha}/(\alpha !)^{1/2}, \alpha \in N^{m}$$ form an orthonormal basis in the Bargmann space. For positive real $\lambda$ the representation $\pi _{\lambda}$ is given by
$$ (\pi_{\lambda}(t,z)F)(w)=F(w-z) exp (i\lambda t+2\lambda(wz-|z|^{2}/2))$$ and for negative $\lambda$ by $\pi_{\lambda}(t,z)=\pi_{|\lambda|}(-t,-z)$. The Fourier transform on $L^{1}(H_{m})$ is given by the formula
$$ \hat{f}(\lambda)=\int_{H^{m}}f(x)\pi_{\lambda}(x)dx, f\in L^{1}(H_{m})$$ and can be extended to an isomorphism between $L^{2}(H_{m})$ and the space of operator valued functions $\hat{f}(\lambda)$ such that
$$\int _{R-0}\|\hat{f}(\lambda)\|^{2}_{HS}\lambda^{k}d\lambda$$ exists. Here $\|.\|_{HS}$ is the Hilbert Schmidt norm.
 If we set $z_{k}=x_{k}+iy_{k}$ then $(t,x_{1},..., x_{m},y_{1},...,y_{m})$ form a coordinate system on $H_{m}$. In this coordinate system we define the following vector fields

$$ X_{k}= \partial_{x_{k}}    +2y_{k}\partial_{t},      1\leq k \leq m $$
$$ X_{k}= \partial_{y_{k}}    -2x_{k}\partial_{t},     m+1\leq k\leq 2m$$
$$T= \partial_{t} .$$
  The fields $T,X_{1},...,X_{2m}$  form the basis for the left-invariant vector fields on $H_{m}$. Every element on $H_{m}$ has a unique representation

$$exp(a_{0}T+a_{1}X_{1}+...+a_{2m}X_{2m}), a_{k}\in R     $$where $exp$ is the exponential map from Lie algebra onto group. One can easely verify that for the fixed integer $j\in Z$ all elements of the form

$$exp(2^{2j}n_{0}T+2^{j}n_{1}X_{1}+...+2^{j}n_{2m}X_{2m})$$ where $n_{k}$ are integers form a discrete subgroup $\Gamma_{j}=\delta_{2^{j}}\Gamma_{0},\Gamma_{0}=\Gamma.$

  The sub-Laplacian $D=-X_{1}^{2}-...-X_{2m}^{2}$ is a second order self-adjoint and positive definite hypoelliptic operator in $L^{2}(H_{m})$  which is homogeneous with respect to the above dilations.

  Using sub-Laplacian $D$ one can introduce the Sobolev scale of spaces with the norm $ \|f\|_{S^{\sigma}(H_{m})}=\|(I+D)^{\sigma /2}f\|, \sigma \geq 0$. As was
shown by Folland [3] (see also [2], [5], [9], [10], [11]) this norm is equivalent to the norm $\|f\|+\|D^{\sigma /2}f\|$ and if $\sigma =r$ is an integer to the norm
$$\|f\|+\sum_{1\leq i_{1},...,i_{r}\leq 2m}\|X_{i_{1}}...X_{i_{r}}f\|.$$
 For negative $\sigma $ spaces $S^{\sigma}(H_{m})$ can be introduced using duality. The full scale $S^{\sigma}(H_{m}), -\infty <\sigma < \infty $ serves the sub-Laplacian $D$ in the same way as standard Sobolev spaces $ H^{\sigma}(R^{d}), -\infty <\sigma <\infty $ serve standard Laplacian $\Delta $.

  {\bf 2.} First of all we introduce an abstract definition of Paley-Wiener functions.

  Let  $E$ be a Hilbert space with the norm $\|.\|$ and $D$ a self-adjoint positive definite operator in $E.$ According to the spectral theory [7] there exist a direct integral of Hilbert spaces $X=\int X(\lambda )dm (\lambda )$ and a unitary operator $F$ from $E$ onto $X$, which transforms domain of $D^{k}$ onto $X^{k}=\{x\in X|\lambda ^{k}x\in X\}$ with norm

$$\|x(\lambda )\|_{X^{k}}= \left (\int^{\infty}_{0} \lambda ^{2k} \|x(\lambda )\|^{2}_{X(\lambda )} dm (\lambda ) \right )^{1/2}  $$
besides $F(D^{k} f)=\lambda ^{k} (Ff), $ if $f$ belongs to the domain of $D^{k}$.
As known, $X$ is the set of all $m $-measurable functions $\lambda \rightarrow x(\lambda )\in X(\lambda ) $, for which the norm

$$\|x\|_{X}=\left ( \int ^{\infty }_{0}\|x(\lambda )\|^{2}_{X(\lambda )} dm (\lambda ) \right)^{1/2} $$
is finite.

We will say that a vector $f$ from $E$ belongs to $PW_{\omega}(D)$ if its "Fourier transform" $Ff$ has support in $[0 , \omega ] $. The next theorem can be considered as an abstract version of Paley-Wiener theorem.

  THEOREM 1. The following conditions are equivalent:

a) a vector $f$ belongs to $PW_{\omega}(D)$;

b) a vector satisfies Bernstein inequality

$$\|D^{k}f\| \leq \omega ^{k}\|f\|  $$

for every natural $k$;

  PROOF. Let $f$ belongs to the space $PW_{\omega }(D)$ and $Ff=x\in X$.
  Then

$$\left(\int ^{\infty}_{0}\lambda ^{2k}\|x(\lambda)\|^{2}_{X(\lambda)}dm
(\lambda)\right)^{1/2}=\left(\int^{\omega}_{0}\lambda^{2k}\|x(\lambda )\|
^{2}_{X(\lambda)}dm(\lambda)\right)^{1/2}\leq \omega ^{k}\|x\|_{X} , k\in N , $$
 which gives Bernstein inequality for $f$.

Conversely, if $f$ satisfies Bernstein inequality then $x=Ff$ satisfies $\|x\|
_{X^{k}}\leq \omega^{k}\|x\|_{X}.$ Suppose that there  exists a set $\sigma \subset [ 0 , \infty ]\setminus [ 0, \omega]$ whose $m $-measure is not zero and $x|_{\sigma }\neq 0.$ We can assume that $\sigma \subset [\omega +\epsilon , \infty )$ for some $\epsilon >0.$ Then for any $k\in N$ we have

$$ \int _{\sigma }\|x(\lambda )\|^{2}_{X(\lambda)}dm (\lambda )\leq \int
^{\infty }_{\omega +\epsilon }\lambda ^{-2k}\|\lambda ^{k}x(\lambda)\|^{2}
_{X(\lambda)}d\mu \leq \|x\|^{2}_{X}\left(\omega /\omega +\epsilon \right)
^{2k}$$ which shows that or $x(\lambda)$ is zero on $\sigma $ or $\sigma $
 has measure zero.

 It is eviden that the set $\bigcup _{
\omega >0}PW_{\omega }(D)$ is dense in $E$ and that the $PW_{\omega }(D)$ is
a linear closed subspace in $E$.

  In the case of a stratified group $H_{m}$ we use sub-Laplacian $D$ in the
   space $L_{2}(H_{m})$. It is a self-adjoint positive definite operator. We
   apply the above construction to the operator $D $ and it gives us the notion of the space $PW_{\omega}(D)$ on the group $H_{m}$.
  Using results from [4] one can show that $f$  belongs to $PW_{\omega}(D)$ if
   and only if its Fourier transform has compact support in the following
   sense:
$${\hat{f}(\lambda)}_{\alpha,\beta}=0,  (2|\beta|+m)|\lambda|>\omega ^{2}$$ where ${\hat{f}(\lambda)}_{\alpha,\beta}=(\hat{f}(\lambda)F_{\alpha,\lambda},F_{\beta,\lambda})$ and inner product is taken in the Bargmann space.

Let $B(x,r)$ be a ball in homogeneous metric $\rho (g,h)=|X-Y|, g=exp X, h=exp Y$ with center $x\in H_{m}$ and radius $r$. Suppose that $\left\{B(x_{\gamma}, r)\right \}_{x_{\gamma}\in \Gamma}$ is a cover of $H_{m}$. It is clear that this cover has a finite multiplicity $M$ in the sense that every ball from this family has non-empty intersections with no more than $M$ other balls from the same family. Since metric $\rho (x,y)$ is homogeneous the family of balls $\left \{B(x_{\gamma}, 2^{j}r)\right \}_{x_{\gamma}\in \Gamma _{j}}$  will  also be a cover of $H_{m}$ of the same multiplicity $M$.

   Given a subgroup $\Gamma _{j}$ and a sequence $\{s_{\gamma}\}\in l_{2}$ we will be interested to find a function $s_{k,j}\in S^{2k}(H_{m}), k >Q/4 $ such that

a)$  s_{k,j}(x_{\gamma})=s_{\gamma}, x_{\gamma}\in \Gamma_{j};$

b) function $s_{k,j}$ minimizes functional $u\rightarrow \|D^{k}u\|$.

  Let us remark first that the same problem for functional $u\rightarrow \|u\|_{S^{2k}(H_{m})}, u\in S^{2k}(H_{m}), k>Q/4$ can be solved easely.

 Pick a ball $B(0,r)$ of very small radius $r$ and then by translations construct the familly of pair ways disjoint balls $B(x_{\gamma},r), x_{\gamma}\in \Gamma_{j}$. In the ball $B(0,r) $ we consider any function $\varphi _{0} \in C^{\infty }_{0}(B(0,r))$ such that $\varphi _{0}(0)=1$. Using translations we construct similar functions $\varphi _{\gamma }$ in balls $B(x_{\gamma },r)$. Because of invariance all this functions have the same Sobolev norm
$$ \|\varphi _{\gamma}\|_{S^{k}(H_{m})}=\|\varphi _{\gamma}\|+\sum _{1\leq i_{1}\leq ... \leq i_{k}\leq 2 m}\|X_{i_{1}}X_{i_{2}}... X_{i_{k}}\varphi _{\gamma}\|, \gamma =1,2, ... .$$

It is clear that for any sequence $\{s_{\gamma}\}\in l_{2}$ the formula
$$f=\sum s_{\gamma}\varphi _{\gamma}$$ defines a function from $S^{k}(H_{m}).$  Let $Pf$ will denote the orthogonal projection of this function $f$ ( in the Hilbert space $ S^{2k}(H_{m})$ with natural inner product) on the subspace $U^{2k}(\Gamma_{j})=\{f\in S^{2k}(H_{m})|f(x_{\gamma})=0\}$ with $ S^{2k}(H_{m})$-norm. Then the function $g=f-Pf$ will be a unique solution of the above minimization problem for the functional $u\rightarrow \|u\|_{ S^{2k}(H_{m})}, k>Q/4$.

  The problem with functional $u\rightarrow \|D^{k}u\|$ is that it is not a norm. But fortunatly we are able to show that  for all natural $k>Q/4$ and all integer $j$ the norm

$$\|D^{k}f\|+\left (\sum_{x_{\gamma}\in \Gamma_{j}}|f(x_{\gamma})|^{2}\right)^{1/2}$$ is equivalent to the norm $\|f\|_{S^{2k}(H_{m})}.$  So, the above procedure can still be applied to the Hilbert space $S^{2k}(H_{m})$ with inner product $$<f,g>=\sum_{x_{\gamma}\in \Gamma_{i}}f(x_{j})g(x_{j})+<D^{k/2}f,D^{k/2}g>$$ and it clearly proves existance and uniqueness of the solution of our minimization problem for the functional $u\rightarrow \|D^{k}u\|, k>Q/4$.\bigskip

 { \bf 3. } We will need the following lemmas.

  LEMMA 2. If $A$ is a self-adjoint operator in a Hilbert space and for some element $f$
$$\|f\|\leq b+ a\|Af\|, a>0,$$
then for all $m=2^{l}, l=0,1,2, ...$
$$\|f\|\leq mb+8^{m-1}a^{m}\|A^{m}f\|$$
as long as $f$ belongs to the domain of $A^{m}$.

  PROOF. For any self-adjoint operator $B$ in a Hilbert space we
have

$$\|f\| \leq \|(I+\varepsilon iB)f\|$$ and the same for the operator $(I-\varepsilon iB)$. It gives

$$\varepsilon \|Bf\| \leq \|(I-\varepsilon iB)f\|+\|f\| \leq \|(I+\varepsilon ^{2}B^{2})f\|+\|f\|\leq \varepsilon^{2}\|B^{2}f\|+2\|f\|.$$

 So, for any $f$ from the domain of $B^{2}$ we have the inequality

 $$ \|Bf\|\leq \varepsilon \|B^{2}f\|+2/\varepsilon \|f\|,   \varepsilon > 0 . $$

Our lemma is true for $m=1$. If it is true for $m$ then applaing the last inequality for $B=A^{m}$ we obtain

$$\|f\|\leq mb + 8^{m-1}a^{m}(\varepsilon\|A^{2m}f\|+2/\varepsilon\|f\|).$$

Setting $\varepsilon=8^{m-1}(a)^{m}2^{2}$, we obtain

$$\|f\| \leq 2mb + 8^{2m-1}(a)^{2m}\|A^{2m}f\|.$$

The lemma 2 is proved.

We consider Sobolev spaces $S^{\sigma}(H_{m})$ with the norm $\|f\|_{S^{\sigma}(H_{m})}=\|f\|+\|D ^{\sigma /2}f\|, \sigma >0$
and for any open $\Omega $ in $H_{m}$ we define the space $S^{\sigma }(\Omega)$ as the collection of all restrictions $g_{\Omega}=g|_{\Omega}, g\in S^{\sigma }(H_{m})$ with the norm $\|g_{\Omega}\|_{S^{\sigma }(\Omega)}=\inf \|g\|_{S^{\sigma }(H_{m})}$ where $g$ runs over the set of all functions from $S^{\sigma }(H_{m})$ whose restriction to $\Omega$ gives $g_{\Omega}$. Let $B(\lambda,M)={B(x_{\gamma}, \lambda)}$ be a cover of $H_{m}$ of finite multiplicity $M$. We introduce a map
$$T_{B(\lambda,M)}:S^{\sigma }(H_{m}) \rightarrow l_{2}(S^{\sigma }(B_{\gamma})), \sigma \geq 0,$$
$$T_{B(\lambda,M)}(g)=\left\{g_{\gamma}\right \}, g_{\gamma}=g|_{B(x_{\gamma},\lambda)}$$
where the Hilbert space on the right is defined as the set of all sequences
${g_{\gamma}}, g_{\gamma}\in S^{\sigma }(B(x_{\gamma},\lambda))$ for which
$(\sum_{\gamma}\|g_{\gamma}\|^{2}_{S^{\sigma }(B(x_{\gamma},\lambda)} )^{1/2}
<\infty$.

  LEMMA 3. For any natural $M$ and any $\sigma \geq 0$ there exists a $C=C(M,\sigma)$ such that for every cover $B(\lambda ,M), \lambda >0$,
 $$\|T_{B(\lambda,M)}\| \leq C(M, \sigma ) max (1,\lambda^{-\sigma }).$$

  PROOF. Let $\theta \in C^{\infty}_{0}(R), \theta (t)=1, |t|\leq 1, supp \>\theta \subset [-2, 2].$ We define $\theta _{\lambda , \gamma }(x)=\theta (\rho (0, \delta_{\lambda ^{-1}}(xx^{-1}_{\gamma})) , x\in H_{m}, \lambda >0.$ Then $\theta _{\lambda }\in C^{\infty}_{0}(H_{m}), \theta _{\lambda}(x)=1, x\in B(x_{\gamma}, \lambda ), supp \>\theta _{\lambda } \subset B(x_{\gamma}, 2\lambda ).$ It is clear that $|X_{i_{1}} ... X_{i_{r}} \theta_{\lambda }(x)|\leq C(r, \theta )\lambda ^{-r}.$ Therefore if $f\in S^{k}(H_{m}), k\geq 0 $ is an integer, then
$$\|f|_{B(x_{\gamma },\lambda)}\|^{2}_{S^{k}(B(x_{\gamma}, \lambda))}\leq\|f\theta _{\lambda}\|^{2}_{S^{k}(H_{m})}\leq $$ $$  \sum _{|r|\leq k} \int _{B(x_{\gamma}, 2\lambda )}|X_{i_{1}} ... X_{i_{r}}(f\theta_{\lambda})(x)|^{2}d\mu \leq $$ $$C(k,\theta )max(1,\lambda ^{-2k}) \sum _{|r|\leq k}\int _{B(x_{\gamma}, 2\lambda )}|X_{i_{1}} ... X_{i_{r}}f(x)|^{2}d\mu  $$

and then
$$ \sum_{\gamma}\|f|_{B(x_{\gamma },\lambda)}\|^{2}_{S^{k}(B(x_{\gamma},\lambda ))}\leq $$ $$C(k,\theta)max(1,\lambda ^{-2k})\sum_{\gamma} \sum_{|r|\leq k}\int _{B(x_{\gamma},2\lambda)}|X_{i_{1}}...X_{i_{r}}f(x)|^{2}d\mu \leq $$ $$C(k,M,\theta)max(1,\lambda^{-2k})\|f\|^{2}_{S^{k}(H_{m})}.$$
Thus for natural $s=k$ lemma is proved. General case can be obtained by interpolation since for the complex interpolation functor $[ . , . ]_{\theta}$
$$\left[l_{2}(L_{2}(B_{\gamma})),l_{2}(S^{\sigma }(B_{\gamma}))\right]_{\theta}=l_{2}(S^{\theta \sigma }(B_{\gamma})), 0<\theta < 1.$$

 The proofs of all main results in the present paper are based on the following inequalities.

  THEOREM 4. There exist a $j_{0}\in Z$ and a constant  $ C_{0} \geq 0 $ such that for $ j\leq j_{0}$ and every $ f \in S^{2k}(H_{m}),k=2^{l}Q, l= 1,2, ... , $ the following inequality takes place
$$\|f\|\leq 2^{l}C_{0}\left (\sum_{x_{\gamma}\in \Gamma _{j}}|f(x_{\gamma})|^{2}\right )^{1/2}+(C_{0}2^{j/2Q})^{k}\|D^{k}f\|.$$

  In particular for $f\in U^{k}(\Gamma_{j})$

$$\|f\|\leq (C_{0}2^{j/2Q})^{k}\|D^{k}f\|.$$

 PROOF. Let $\{B(x_{\gamma},1)\}_{x_{\gamma}\in\Gamma}$ be a cover of $H_{m}$ of the finite multiplicity $M$. The cover   $\{B(x_{\gamma},2^{j})\}_{x_{\gamma}\in\Gamma_{j}}, \Gamma _{j}=\delta_{2^{j}}\Gamma $ also has the same multiplicity $M$. Let ${\psi _{\gamma}}, supp \psi _{\gamma} \subset B_{\gamma}$ be a corresponding partition of unity.

 For a function $f$ from $ S^{\sigma}(H_{m}), \sigma >Q/2$ we consider decomposition
$$f(x)=\sum _{\gamma} f(x)\psi_{\gamma}(x)=\sum _{\gamma} f(x_{\gamma})\psi_{\gamma}(x)+\sum _{\gamma}(f(x)-f(x_{\gamma}))\psi_{\gamma}(x)$$ and then
$$\|f\|^{2}\leq C\left \{ \sum _{\gamma}|f(x_{\gamma})|^{2}+\sum _{\gamma}
\int_{B(x_{\gamma},\lambda)}|f(x)-f(x_{\gamma})|^{2}d\mu \right \},$$ where $C$ depends only on multiplicity $M$.

  Since every vector field on the group $H_{m}$ is a linear combination over $C^{\infty}$ of the fields $ X_{i},[X_{i_{1}}, X_{i_{2}}], 1\leq i,i_{1},i_{2}\leq 2 m$, the Newton-Leibnitz formula gives
$$|f(x)-f(x_{\gamma})|^{2}\leq$$ $$ C 4^{j}\left( \sum _{1\leq i_{1},i_{2}\leq 2m}(\sup_{y\in B(x_{\gamma},2^{j})}|X_{i_{1}}X_{i_{2}}f(y)|)^{2}+\sum_{1\leq i\leq 2m}(\sup_{y\in B(x_{\gamma},2^{j})}|X_{i}f(y)|)^{2} \right).$$

Applying anisotropic version of Sobolev inequality [3] we obtain
$$|f(x)-f(x_{\gamma})|^{2}\leq $$ $$C4^{j}\left (\sum_{1\leq i_{1},i_{2}\leq 2m}(\sup_{y\in B(x_{\gamma},2^{j})}|X_{i_{1}}X_{i_{2}}f(y)|)^{2}+  \sum_{1\leq i\leq 2m}(\sup_{y\in B(x_{\gamma},2^{j})}|X_{i}f(y)|)^{2} \right )\leq$$

 $$C4 ^{j} \left ( \sum _{1\leq i_{1},i_{2}\leq 2m}\| X_{i_{1}}X_{i_{2}}f\|^{2}_{S^{Q/2+\varepsilon}(B(x_{\gamma},2^{j}))}+\sum_{1\leq i\leq 2m}\|X_{i}f\|_{S^{Q/2+\varepsilon}(B(x_{\gamma} ,2^{j})}^{2} \right),$$ where $ x\in B(x_{\gamma}, 2^{j}),\varepsilon > 0, C=C(X_{i},...,X_{2m};\varepsilon).$

An application  of lemma 3 gives
$$\sum_{\gamma}\int _{B(x_{\gamma},2^{j})}|f(x)-f(x_{\gamma})|^{2}d\mu \leq $$ $$C (2^{j})^{Q+2} \left (\sum _{1\leq i_{1},i_{2}\leq 2m}\sum_{\gamma}\|X_{i_{1}}X_{i_{2}}f\|^{2}_{S^{Q/2+\varepsilon}(B(x_{\gamma},2^{j}))}+  \sum_{1\leq i\leq 2m}\sum _{\gamma}\|X_{i}f\|^{2}_{S^{Q/2+\varepsilon}(B(x_{\gamma} ,2^{j})} \right) \leq$$

 $$ C(2^{j})^{2-2\varepsilon}\left ( \sum _{1\leq i_{1},i_{2}\leq m}\|X_{i_{1}}X_{i_{2}}f\|^{2}_{S^{Q/2+\varepsilon}(H_{m}}+  \sum_{1\leq i\leq 2m}\|X_{i}f\|_{S^{Q/2+\varepsilon}(H_{m})}^{2} \right )\leq$$

 $$C(2^{j})^{2-2\varepsilon}\left( \sum _{1\leq i_{1},i_{2}\leq m}\|X_{i_{1}}X_{i_{2}}f\|^{2}_{S^{\sigma}(H_{m})}+ \sum_{1\leq i\leq 2m}\|X_{i}f\|_{S^{\sigma}(H_{m})}^{2} \right),$$ where $\sigma \geq Q/2+\varepsilon$, $C$ depends only on $X_{i},..., X_{2m}$ on $ \sigma$ and on multiplicity $M$. Since
$$\| X_{i}f\|_{S^{\sigma}(H_{m})}+\|X_{i_{1}}X_{i_{2}}f\|_{S^{\sigma}(H_{m})}\leq C\left\{\|f\|+\|D^{1+\sigma/2}f\| \right \},$$ we have for particular choice of $\varepsilon=1/2, \sigma =2Q-2, $
$$\|f\|\leq C \left \{\left(\sum_{\gamma}|f(x_{\gamma})|^{2}\right)^{1/2}+2^{j/2}\|D^{Q}f\|+2 ^{j/2}\|f\|\right \}, $$ where $C$ depends only on $X_{1},...X_{2m}$ and multiplicity $M$. Thus, if $j $ is smaller than some $j_{0}=j_{0}(X_{1},...,X_{2m};M)$ it gives
$$\|f\|\leq C\left \{\left(\sum _{\gamma}|f(x_{\gamma})|^{2}\right)^{1/2}+C2^{j/2}\|D^{Q}f\|\right \}, C=C(X_{1}, ...,X_{2m};M).$$

Using lemma 2 for $A=D^{Q}, b=\left (\sum_{\gamma}|f(x_{\gamma})|^{2}\right)^{1/2}$ we obtain $$\|f\|\leq C2^{l}\left (\sum_{\gamma}|f(x_{\gamma})|^{2}\right )^{1/2}+(C_{0}2 ^{j/2})^{2^{l}}\|D^{2^{l}Q}f\|,$$ where $ l=0,1,2, ...,C_{0}=8C.$ After all, for $f\in U^{2^{l}Q}(\Gamma _{j}), j<j_{0} $ $$\|f\|\leq (C_{0}2 ^{j/2})^{2^{l}}\|D^{2^{l}Q}f\|, l=0,1,2,...$$ Theorem 4 is proved.

  LEMMA 5. For any natural $k>Q/4$ and any $\Gamma_{j}=\delta_{2^{j}}\Gamma, j\in Z,$ the norm $\|f\|_{S^{2k}(G)}$ is equivalent to the norm
$$\|D^{k}f\|+\left(\sum _{x_{\gamma
}\in \Gamma_{j}}|f(x_{\gamma })|^{2}\right )^{1/2}.$$

  PROOF. The proof of the theorem 4 shows that for every natural $k>Q/4$ there exists a $j(k)$ such that for every $j\leq j(k)$ there is a $C=C(k,j)$ for which$$\|f\|\leq C\left \{
\|D^{k}f\|+\left(\sum _{x_{\gamma
}\in \Gamma_{j}}|f(x_{\gamma })|^{2}\right )^{1/2}\right \}.$$

Now using homogeneity arguments one can easely show that for every natural $k>Q/4$ and every integer $j$ there exists a $C=C(k,j)$ for which the above inequality takes place.

In order to prove inverse inequality we consider $C^{\infty }_{0}(H_{m})$ functions $\phi_{\gamma}$ with disjoint supports such that  $\phi_{\gamma}(x_{\gamma})=1$. Using Sobolev embedding theorem we obtain for $k>Q/4$
$$\left(\sum _{\gamma}|f(x_{\gamma})|^{2}\right)^{1/2}\leq C_{k}\left(\sum _{\gamma}\|f\phi _{\gamma}\|^{2}_{S^{2k}(H_{m})}\right)^{1/2}\leq C_{k}\|f\|_{S^{2k}(H_{m})}, k> Q/4.$$
 The proof of the lemma 5 is finished.

  Next, we consider the following minimization problem [6]. For the given $f\in S^{2k}(G), k>Q $, the $s_{k,j}(f)\in S^{2k}(G)$ will be the function that minimizes $u\rightarrow \|D^{k}u\|$ and takes the same values on $\Gamma_{j}$ i.e. $s_{k,j}(f)|_{\Gamma_{j}}=f_{\Gamma_{j}}.$  Since $D$ is invariant with respect to translations it is clear that $s_{k,j}(f)=\sum_{x_{\gamma}\in\Gamma _{j}}f(x_{\gamma})L_{k,j}(xx_{\gamma}^{-1})$ where $L_{k,j}(x)\in S^{2k}(G)$ is the function that minimizes the same functional and $L_{k,j}(0)=1$, and is zero at all other points of $\Gamma _{j}$. In classical case such functions are called Lagrangian splines .\bigskip

 {\bf 4.} We prove the following approximation theorem.

  THEOREM 6. There exists $  c_{0}> 0 $ such that for $j\leq j_{0}$  the following estimate takes place
$$ \|f-s_{k,j}(f)\|\leq (c_{0}2 ^{j/2Q})^{k}\|D ^{k}f\| ,  f\in S^{2k}(H_{m}), k=2^{l}Q, l=1,2, ... .$$

  PROOF. If $f\in S^{2k}(H_{m}), k=2^{l}Q$  then $f-s_{k,j}(f)\in U^{2k}(\Gamma_{j})$ and according to the theorem 4 we have

  $$\|f-s_{k,j}(f)\|\leq (C_{0}2 ^{j/2Q})^{k}\|D^{k}(f-s_{k,j}(f))\|.$$

 Using minimization property of $s_{k,j}(f)$ we obtain

$$\|f-s_{k,j}(f)\|\leq (c_{0}2 ^{j/2Q})^{k}\|D ^{k}f\|, k=2^{l}Q,$$ where $c_{0}=2C_{0}$ and the constant $C_{0}$ is from theorem 4.

 Using theorem 6 and Bernstein inequality we immediatly come to the following uniqueness and reconstruction theorem.

  THEOREM 7. For the same  constant  $c_{0}> 0$ as above

a) every  function $f\in PW_{\omega}(D),  \omega >0$ is uniquely determined by its values on any set $\Gamma_{j}=\delta _{2^{j}}\Gamma$ as long as $j<-2Qlog_{2}(c_{0}\omega );$

b) for every such set $\Gamma_{j}=\delta _{2^{j}}\Gamma$ the sequence of splines
  $$s_{k,j}(f)(x)=\sum _{x_{\gamma}\in \Gamma_{j}} f(x_{\gamma})L_{k,j}(xx_
  {\gamma}^{-1}) , k=2^{l}Q, l=1,2, ... ,$$ converges to $f\in PW_{\omega}(D)$
   in $L^{2}(H_{m})$-norm.\bigskip

 {\bf 5.}  As an concluding remark we will show that functions $s_{k,j}$ have
 the following remarkable property (see [7]).

$$D ^{2k}s_{k,j}=\sum_{x_{\gamma}\in \Gamma _{j}}\alpha _{\gamma }\delta (x_{\gamma }),$$ where $\delta (x) $ is the Dirac measure and $\{ \alpha _{\gamma} \} \in l_{2}$.

Indeed, suppose that $s_{k,j}\in S^{2k}(H_{m})$ is a solution to the minimization problem and $h\in U^{2k}(\Gamma_{j}).$ Then

$$\|D ^{k}(s_{k,j}+\lambda h)\|^{2}=\|D ^{k}s_{k,j}\|^{2}_{2}+2Re {\lambda\int_{H_{m}}D ^{k}s_{k,j} D ^{k}h}d\mu+|\lambda |^{2}\|D ^{k}h\|^{2}_{2}.$$

 The function $s_{k,j}$ can be a minimizer only if for any $h\in U^{2k}(\Gamma_{j})$

$$\int_{H_{m}} D^{k}s_{k,j}D^{k}hd\mu=0.$$

So, the function $\Phi= D ^{k}s_{k,j}\in L_{2}(H_{m}) $ is orthogonal to $ D ^{k}U^{2k}(\Gamma_{j}) $.
Let $\varphi_{\gamma}$ be the same set of functions as above and   $h\in C_{0}^{\infty}(H_{m})$. Then the function $h-\sum h(x_{\gamma})\varphi_{\gamma}$ belongs to the
$U^{2k}(\Gamma_{j})\cap C_{0}^{\infty}(H_{m}).$ Thus,

$$0=\int_{H_{m}}\Phi \overline{D ^{k}(h-\sum h_{\gamma}\varphi_{\gamma})}d\mu=\int_{H_{m}} \Phi \overline{D ^{k} h}d\mu-\sum \overline{h(x_{\gamma})}\int_{H_{m}}\Phi \overline{D ^{k}\varphi_{\gamma}}d\mu.$$

In other words

$$D ^{k}\Phi =\sum_{x_{\gamma}\in \Gamma _{j}}\alpha _{\gamma }\delta (x_{\gamma }),$$ or
$$D ^{2k}s_{k,j}=\sum_{x_{\gamma}\in \Gamma _{j}}\alpha _{\gamma }\delta (x_{\gamma }),$$ where $\delta (x) $ is the Dirac measure.

Moreover for any integer $r>0$

$$\sum _{\gamma =1}^{r}|\alpha _{\gamma }|^{2}=<\sum_{1}^{\infty}\alpha _{\gamma }\delta (x_{\gamma }), \sum _{1}^{r}\alpha _{\gamma} \phi _{\gamma}>\leq C \|\sum _{1}^{\infty} \alpha _{\gamma} \delta (x_{\gamma})\|_{S^{-2k}(H_{m})}(\sum_{1}^{r} |\alpha _{\gamma}|^{2})^{1/2},$$
where $C$ is independent on $r$. It shows that the sequence $\{\alpha _{\gamma }\}$ belongs to $l_{2}$.

$$ACKNOWLEDGMENTS$$
   I thank Professors L. Ehrenpreis and R. Strichartz for helpful conversations. \bigskip

$$References $$

1. J. Benedetto, Irregular sampling and frames, in "Wavelets: A tutorial in Theory and Applications" (C.K.Chui, Ed.), 445-507. Academic Press, Boston,
1992.

2. A.F.M. Ter Elst, D.W. Robinson, Subelliptic operators on Lie groups: regularity, J. Austral. Math. Soc. (Series A) 57 (1994), 179-229.

3. G. Folland, Subelliptic estimates and function spaces on nilpotent Lie groups, Ark. Mat., t. 13, (1975), 161-207.

4. D. Geller, Fourier Analysis on Heisenberg group, Proc. Nat. Acad. Sci. U.S.A. 74(1977)1328-1331.

5. S. Krein, I. Pesenson, Interpolation Spaces and Approximation on Lie Groups, pp. 200, The Voronezh State University, Voronezh, 1990,(in Russian).

6. P. Lemarie, Bases d'ondelettes sur les groupes stratifies,
Bull. Soc. Math. France, 177 (1990), 211-232.

7. J.-L. Lions, E. Magenes, Non-Homogeneous Boundary Value Problem and Applications, Springer-Verlag, 1975.

8. Y. Meyer, Ondelettes et Operateurs, in two volumes, Hermann, Paris, 1990.

9.  I. Pesenson, The Nikol`skii-Besov Spaces in Representations of Lie Groups. Dokl. Acad. Nauk, USSR, v. 273, No. 1 (1983), 45-49; Engl. transl. in Soviet  Math. Dokl. 28 (1983).

10. I. Pesenson, The Bernstein Inequality in the Space of Representation of Lie group. Dokl. Acad. Nauk USSR v. 313, No. 4 (1990), 86-90; Engl. transl. in Soviet Math. Dokl. 42 (1991).

11. I. Pesenson, Approximation in Representation Space of a Lie Group. Izvestiya VUZ, Mathematika v. 34, No.7 (1990), 43-50; Engl. transl. in Soviet Mathematics v. 34, No 7  (1991).

12. I. Pesenson, Lagrangian splines, spectral entire functions and Shannon-Whittaker theorem on manifolds, Temple University Research Report 95-87 (1995),1-28.

13. I. Pesenson, A sampling theorem on homogeneous manifolds, Preprint (1996).

14. I. Schoenberg, "Cardinal Spline Interpolation." CBMS, Vol. 12, SIAM,
Philadelphia, 1973.

\end{document}